\documentclass{amsart}

\usepackage{xspace,amssymb,amsfonts,euscript,eufrak,amsthm,amsmath}
\usepackage{graphicx,ifpdf}
\ifpdf \usepackage{epstopdf}
 \usepackage{pdfsync}\fi
\usepackage{palatino}
\usepackage{color}

 \title[Characteristic cycles and decomposition numbers]
 {Characteristic cycles and decomposition numbers}

\input xy
\xyoption {all}

  \newcommand{\nc}{\newcommand}
  \newcommand{\renc}{\renewcommand} 

\usepackage[latin1]{inputenc}

\nc{\E}{\mathbb{E}}

\nc{\Q}{\mathbb{Q}}

\nc{\triright}{\stackrel{[1]}{\to}}

\def\to{\rightarrow}
\def\from{\leftarrow}
\def\longto{\longrightarrow}

\nc{\Br}{\mathcal{B}}
\nc{\id}{id}
\nc{\HotRR}{{}_R\mathcal{K}_R}
\nc{\HotR}{\mathcal{K}_R}
\nc{\excise}[1]{}
\nc{\defect}{\text{df}}
\nc{\h}[1]{\underline{H}_{#1}}
\nc{\Z}{\mathbb{Z}}
\nc{\R}{\mathbb{R}}
\nc{\C}{\mathbb{C}}
\renc{\P}{\mathbb{P}}
\renc{\O}{\mathcal{O}}
\nc{\N}{\mathbb{N}}

\nc{\F}{\mathcal{F}}
\nc{\G}{\mathcal{G}}

\nc{\nilp}{\mathcal{N}}

\nc{\Ga}{\mathbb{G}_a} 
\nc{\Gm}{\mathbb{G}_m} 

\nc{\Loc}{\mathcal{L}}

\nc{\A}{\mathbb{A}} 

\nc{\D}{\mathbb{D}}
\nc{\ic}{\mathbf{IC}}
\nc{\Ic}{\mathcal{I}\mathcal{C}}

  \newtheorem{thm}{Theorem}[section]
  \newtheorem{lem}[thm]{Lemma}

  \newtheorem{cor}[thm]{Corollary}

  \theoremstyle{remark}
  \newtheorem{remark}{Remark}

\DeclareMathOperator{\CC}{CC}
\DeclareMathOperator{\re}{Re}
\DeclareMathOperator{\rk}{rk}

    \def\CM{{\mathbb{C}}}

    \def\FM{{\mathbb{F}}}

    \def\RM{{\mathbb{R}}}

    \def\ZM{{\mathbb{Z}}}

    \def\DC{{\mathcal{D}}}

    \def\LC{{\mathcal{L}}}
    \def\MC{{\mathcal{M}}}
    
    \def\OC{{\mathcal{O}}}

\def\FS{{\EuScript F}}

\def\LS{{\EuScript L}}
\def\MS{{\EuScript M}}

\def\G{\Gamma}

\def\D{\Delta}

\def\l{\lambda}
\def\L{\Lambda}
\def\O{\Omega}

\newcommand{\Mod}{\textrm{-Mod}}

\nc{\simto}{\stackrel{\sim}{\to}}
\nc{\Perv}{{\mathbf{P}}}
\def\gr{\mathrm{gr}}

\nc{\cmtk}[1]{\color{red}{{\fbox{K}} #1}\color{black}}

\begin{document}

\author{Kari Vilonen}
\address{Department of Mathematics, Northwestern University, Evanston,
  IL 60208, USA
  \\
and \\ Department of Mathematics, Helsinki University, Helsinki, Finland}
\thanks{K.\ Vilonen was supported by NSF and by DARPA via AFOSR grant
  FA9550-08-1-0315.\\
Both authors were supported by the EPSRC and would like to thank the
MPI, Bonn for a productive research environemnt.}
\email {vilonen@math.northwestern.edu, kari.vilonen@helsinki.fi}
\author{Geordie Williamson}
\address{Max-Planck-Institut f\"ur Mathematik, Vivatsgasse 7, 53111 Bonn, Germany}
\email{geordie@mpim-bonn.mpg.de}

\maketitle

\section{Introduction}

In the articles  \cite{KL} and \cite{L} Kazhdan and Lusztig made
conjectures that amount to the irreducibility of characteristic
varieties of intersection cohomology sheaves of Schubert varieties in
type $A_n$. In \cite{KS} Kashiwara and Saito provide a counterexample
to this conjecture for $n=7$. Previously  Kasiwara and Tanisaki \cite{KT,T} gave
examples of reducible characteristic varieties in other types. On the other hand, in the context modular
representation theory, Braden discovered examples of non-trivial
decomposition numbers for perverse sheaves (with coefficients in a
finite field) on flag varieties. One of these examples occurs in type
$A_7$. The second author \cite{W2} pointed out recently that similar considerations
give a counterexample to the Kleshchev-Ram conjecture \cite{KR}.

The
purpose of this paper is to observe that this is no accident. Given a
sheaf with integer coefficients extension of scalars yields a sheaf with
rational coefficients and modular reduction yields a sheaf with coefficients in any finite
field. It is an easy observation that all the sheaves in question have
the same characteristic variety. Hence, if a sheaf becomes reducible
modulo a prime $p$ then its characteristic variety is
reducible. Similarly, if a sheaf has irreducible
characteristic variety then it remains irreducible modulo all primes.

In section 2 we state the main result. The proof is given in section
3. We include a rather lengthy discussion of characteristic cycles and
their basic properties from various points of view which we hope will
help clarify the notion. The most natural context for characteristic
cycles is the setting of real (analytic) manifolds and real
constructible sheaves which we discuss briefly in section~\ref{real
  case}. In section 4 we discuss some examples.  

In 2010 on the Isle of Skye the second author gave a talk about
the first examples (in types $B_2$ and $A_7$) of decomposition numbers of
intersection cohomology complexes of Schubert varieties. Michael
Finkelberg remarked that these examples coincided with examples due to
Tanisaki and Kashiwara-Saito of reducible characteristic varieties,
and asked if there was a relation. We would like to thank him 
for his observation and question. We would also like to thank Tom Braden for useful
correspondence and verifying that a certain Schubert variety in type
$A_7$ has reducible characteristic variety, which led us to
believe beyond reasonable doubt that there was something going on!

\section{Statement of the theorem} \label{sec:statement}

Let $X$ denote a smooth complex algebraic variety equipped with a Whitney
stratification
\[
X = \bigsqcup_{\l \in \L} X_\l
\]
by locally closed connected smooth subvarieties. To each stratum $X_\l
\subset X$ and local system $\LS$ of free 
$\Z$-modules on $X_\l$ one may associate two objects:
\begin{enumerate}
\item[i)] $\Ic(\overline{X_\l}, E_\LS)$: the intersection cohomology
  $\DC_X$-module, which is a regular holonomic $\DC_X$-module extending
  $E_\LS$, a vector bundle with flat connection on $X_\l$ whose monodromy
  is given by $\LS \otimes_\ZM \CM$.
\item[ii)] $\ic(\overline{X_\l}, \LS)$: the extension by zero to $X$ of the
  intersection cohomlogy complex extending the local system  $\LS$ on
  $X_{\l}$. This is an object in $\Perv_{\Lambda}(X, \Z)$, the
  abelian category of $\L$-constructible perverse sheaves on $X$ with
  coefficients in $\Z$.
\end{enumerate}

Associated to the above two objects are two basic problems:
\begin{enumerate}
\item[i)] On $\Ic(\overline{X_\l},E_\LS)$ one may choose a good filtration
  $F$. The associated graded is then a sheaf of modules over $\gr \DC_X =
  \OC_{T^*X}$ with support
  contained in the union of the conormal bundles to the strata.
  Taking multiplicities along each conormal bundle gives rise to the
  characteristic cycle
\[
\CC(\Ic(\overline{X_\l},E_\LS)) = \sum_\mu m^{\LS}_{\l\mu} [\overline{T_{X_\mu}^* X}].
\]
A basic problem in the theory of $\DC$-modules is the calculation of the
multiplicities $m^{\LS}_{\l\mu}$.

\item[ii)] Fix a prime number $p$ and consider $\ic(\overline{X_\l}, \LS) \stackrel{L}{\otimes}_\Z \FM_p$
  the reduction modulo $p$ of $\ic(\overline{X_\l}, \LS)$. It is a
  perverse sheaf with $\FM_p$-coefficients on $X$ which in general
  will not be simple. Given any local system $\MS$ of
  $\FM_p$-vector spaces on $X_\mu$ the second
  basic problem asks for the decomposition number:
\[
d^p_{\l,\LS,\mu,\MS} = [ \ic(\overline{X_\l}, \LS) \stackrel{L}{\otimes}_\Z \FM_p : \ic(\overline{X_\mu}, \MS)].
\]
This question is related to interesting questions in representation
theory. For example, in both the geometric Satake equivalence \cite{MV} and the
modular Springer correspondence \cite{J} the above (topological) problem is equivalent to the
(algebraic) problem of determining the decomposition numbers of standard modules.
\end{enumerate}

Let $\rk \MS$ denote the rank of the local system $\MS$. 
Then:

\begin{thm} \label{thm:main}
  $\rk \MS \cdot d^p_{\l, \LS, \mu, \MS} \le m^{\LS}_{\l\mu}$.
\end{thm}

Recall that the characteristic variety of $\Ic(\overline{X_\l},E_\LS)$ is
defined as the support of $\CC(\Ic(\overline{X_\l}),E_\LS)$. The above
theorem has the following consequence:

\begin{cor}
  Suppose that the characteristic variety of $\Ic(\overline{X_\l},E_\LS)$ is
  irreducible. Then $\ic(\overline{X_\l}, \LS)
  \stackrel{L}{\otimes}_\Z \FM_p \cong \ic(\overline{X_\l}, \LS
   \stackrel{L}{\otimes}_\Z \FM_p)$ for all primes $p$. Hence $\ic(\overline{X_\l}, \LS)
  \stackrel{L}{\otimes}_\Z \FM_p$ is simple if $ \LS
   \stackrel{L}{\otimes}_\Z \FM_p$ is.
\end{cor}

\section{Characteristic cycles and the proof}

\subsection{The Euler characteristic:} Let $\Bbbk$ denote either $\Z, \CM$ or $\FM_p$ and let $D^b_c(pt, \Bbbk)$
denote the bounded derived category of constructible sheaves of
$\Bbbk$-modules on a point, also known as the full subcategory of the
bounded derived category of $\Bbbk$-modules consisting of complexes with finitely
generated cohomology. Let $K_\Bbbk$
 denote the Grothendieck group of $D^b_c(pt, \Bbbk)$. We have a canonical
 isomorphism
\[
\chi : K_\Bbbk\simto \Z
\]
fixed by declaring $\chi([\Bbbk]) = 1$, where $[\Bbbk]$ denote the class of
the free module of rank 1 (placed in degree zero). Of course $\chi$ is just the Euler
characteristic if $\Bbbk$ is a field, and is the alternating sum of the
ranks if $\Bbbk = \Z$.

Consider the functors of extension of scalars:
\[
D^b_c(pt, \FM_p)
\stackrel{(-)\stackrel{L}{\otimes}_\Z\FM_p}{\longleftarrow}
D^b_c(pt, \Z) \stackrel{{(-) \otimes_{\Z} \CM}}{\longrightarrow} D^b_c(pt, \CM).
\]
These functors are exact (in the triangulated sense) and hence induce maps on the corresponding Grothendieck groups:
\[
K_{\FM_p} \stackrel{f}{\longleftarrow} K_\Z \stackrel{c}{\longrightarrow} K_{\CM}
\]
These functors preserve the class of free module of rank 1. Hence:

\begin{lem} \label{lem:commute}
 $f$ and $c$ are isomorphisms commuting with $\chi$.
\end{lem}

\subsection{The characteristic cycle of a constructible sheaf:}
Recall that $X = \bigsqcup
X_\l$ is a smooth Whitney stratified complex algebraic variety.
We let $T^*X$ denote the cotangent bundle of $X$, $\pi : T^*X \to X$ the canonical
projection and $T_\l^*X := T_{X_\l}^*X$ the conormal bundle to $X_\l
\subset X$.
A covector $\xi \in T_\l^*X$ is called non-degenerate if $\xi \notin
\overline{T_\mu^* X}$ for all $\mu \ne \l$.

As above let $\Bbbk \in \{ \ZM, \CM, \FM_p \}$ and denote by
$D^b_\Lambda(X;\Bbbk)$ the $\Lambda$-constructible derived category: this is
the full subcategory of the bounded derived category of sheaves of
$\Bbbk$-modules consisting of complexes whose cohomology sheaves
are constructible with respect to the stratification $\Lambda$.

Now fix a constructible complex $\FS \in D^b_\L(X;\Bbbk)$. To $\FS$ we
wish to associate its characteristic cycle, which is a $\Z$-linear
combination of the cycles $[\overline{T_\l^*X}]$ for all $\l \in
\L$. To this end fix a stratum $X_\l$, a point $x \in X_\l$ and a non-degenerate
covector $\xi \in T_\l^*X$ such that $\pi(\xi) = x$. Then there exists
a neighbourhood $U$ of $x$ and a holomorphic function $\phi : U \to
\CM$ such that
\begin{enumerate}
\item[i)] $\phi(x) = 0$ and $d\phi_x = \xi$,
\item[ii)] The image of $d\phi$ intersects $ \bigcup T_\mu ^*X$ only at $\xi$,
  and this intersection is transverse.
\end{enumerate}
Now set
\[
A_{\xi}(\FS) := R\G_{\{\re \phi \ge 0\}}(\FS)_x  \in D^b_\L(\{x\};\Bbbk)
\]
where $\{ \re \phi \ge 0 \} := \{ x \in X \; | \; \re \phi(x) \ge 0
\}$. The complex $A_\xi(\FS)$ is independent of the choice
of $\phi$ and its Euler characteristic
\[
m_\l(\FS) := \chi(A_\xi(\FS))
\]depends
only $\l$. One then defines the characteristic cycle as
\[
\CC(\FS) := \sum_{\l \in \L} m_\l(\FS) [\overline{T_\l^*X}] \in \bigoplus_{\l\in\L}
\ZM [\overline{T_\l^*X}].
\]
For a more concrete description of $A_{\xi}(\FS)$, see for example \cite[\S 2]{SV1}.

\subsection{Basic properties of the characteristic cycle:}
\begin{enumerate}
\item[i)] The characteristic cycle factors through the Grothendieck group of $D^b_\L(X; \Bbbk)$:
  given a distinguished triangle
\[
\FS ' \to \FS \to \FS'' \stackrel{[1]}{\longto}
\]
then $\CC(\FS) = \CC(\FS') + \CC(\FS'')$. (This is immediate from the
definition).
\item[ii)] Given $\FS \in D^b_\L(X, \Bbbk)$ and $X_\l$ open in the
  support of $\FS$ then for any non-degenerate $\xi \in T^*_\l X$ we have
\[
A_\xi(\FS) \cong \FS_x[-d_\l]
\]
where $d_\l$ denotes the complex dimension of $X_\l$. (This may be
checked directly from the definition.) In particular, if $\LS$
denotes a local system of free $\Bbbk$-modules on $X_\l$ and 
$\ic(\overline{X_\l};\LS)$ denotes the intersection cohomology
extension then
$m_\l(\ic(\overline{X_\l};\LS)) = \rk \LS$, the rank of $\LS$.

 \item[iii)] Let $\Perv_\L(X;\Bbbk) \subset D^b_\L(X;\Bbbk)$ denote the abelian
  subcategory of perverse sheaves (see \cite{BBD}). Then
  $\FS \in D^b_\L(X;\Bbbk)$ is perverse if and only if for all $\l \in \L$
  and non-degenerate $\xi \in T^*_\l X$ the complex $A_\xi(\FS)$
is concentrated in degree 0, see \cite[Theorem 10.3.12]{KSh}. 
Indeed, one can argue as in \cite[Theorem 6.4]{GM} to show that $A_\xi(\FS)$ is concentrated in
degree 0 if $\FS$ is perverse. The other direction follows from the fact that if $A_\xi(\FS)=0$ for all  non-degenerate $\xi \in T^*_\l X$ then $\FS=0$.  In particular, for $\FS
\in \Perv_\L(X;\Bbbk)$ we have
\[
\CC(\FS) \in \bigoplus_\l \ZM_{\ge 0} [ \overline{T_\l^*X}].\]

\item[iv)] The characteristic cycle operation commutes with the Riemann-Hilbert
  correspondence. Let $\DC_X\Mod_{\L}$ denote the abelian category of regular
  holonomic $\DC_X$-modules whose characteristic variety is contained
  in $\bigcup \overline{T_\l^*X}$. Recall that the Riemman-Hilbert
  correspondence gives an equivalence of abelian categories
\[
RH : \DC_X\Mod_{\L} \simto \Perv_\L(X;\CM).
\]
Given any $\DC_X$-module $\MC$ one can define its characteristic
cycle using the theory of good filtrations. We have
\[
\CC(\MC) = \CC(RH(\MC)).
\]
This is a theorem due to Kashiwara \cite[\S 8.2]{K}. See the final pages of
\cite{SV2} for a short proof of this fact. The basic idea is that
using results of Ginzburg \cite{G} and \cite{SV1}
both sides can be seen to agree if $\MC$ is the $\DC$-module direct
image of a vector bundle with flat connection on a locally closed
subvariety, and such direct images generate the Grothendieck group of
all holonomic $\DC$-modules.
\item[v)]  The characteristic cycle commutes with extension of scalars:
If $\Bbbk \in \{
\CM, \FM_p\}$ and $\FS \in D^b_\L(X;\ZM)$ then
\begin{equation*} \label{eq:CCreduc}
\CC(\FS) = \CC(\FS \otimes^L_\ZM \Bbbk). 
\end{equation*}
This is immediate from Lemma \ref{lem:commute}
and the fact that the functors $A_\xi(-)$ 
and $(-) \otimes^L_\ZM \Bbbk$ commute up to natural isomorphism.
\end{enumerate}

\subsection{Proof of the theorem:} Fix $\l \in \L$ and a local system
$\LS$ of free $\ZM$-modules on $X_\l$. Let 
$\ic(\overline{X_\l}, \LS)$ denote the intersection cohomology complex
of $\overline{X_\l}$ with coefficients in $\LS$. We have
\[
\ic(\overline{X_\l}, \LS) \otimes_{\ZM} \CM \cong \ic(\overline{X_\l},
\LS \otimes_\ZM \CM)
\]
Fix a prime $p$. In the Grothendieck group of $\Perv_\L(X,\FM_p)$ we
can write
\[
[\ic(\overline{X_\l}, \LS) \otimes_\ZM^L \FM_p  ] = \sum d^p_{\l, \LS,
  \mu, \MS}
[\ic(\overline{X_\mu};\MS)] 
\]
where the sum runs over all pairs $\mu, \MS$ where $\mu \in \L$ and
$\MS$ is a irreducible local system of $\FM_p$-vector spaces on $X_\mu$.
Applying $\CC$ and using i), ii) and iii) above we are led to the inequalities
\[
m_\mu(\ic(\overline{X_\l}, \LS) \otimes_\ZM^L \FM_p) \ge \rk \MS \cdot
d^p_{\l, \LS,
  \mu,\MS}.
\]
By iv) and v) we get (with notation as in \S \ref{sec:statement})

\begin{multline*}
\CC(\Ic(\overline{X_\l},E_\LS)) = \CC(\ic(\overline{X_\l}; \LS
\otimes_\ZM \CM) =
\\
\CC(\ic(\overline{X_\l}, \LS) \otimes \CM)) = \CC(\ic(\overline{X_\l};\LS) \otimes_\ZM^L \FM_p).
\end{multline*}

Hence
\[
m^{\LS}_{\l \mu} = m_\mu(\ic(\overline{X_\l}, \LS) \otimes_\ZM^L \FM_p)
\]
which completes the proof.

\qed

\subsection{Remarks on the real analytic case}\label{real case} Let
$X$ to be a real analytic manifold and consider $D^b_c(X;\Bbbk)$ the
bounded derived category of sub-analytically constructible 
sheaves. In this case the characteristic cycle has been defined by
Kashiwara \cite[\S 3]{K}.
It coincides with the previous definition in the complex case.

In the real case we get an isomorphism
\[
\CC : K(X;\Bbbk) \simto \LC_X
\]
where $K(X;\Bbbk)$ denotes the Grothendieck group of
$D^b_c(X;\Bbbk)$ and $\LC_X$ denotes the group of $\RM^+$-invariant, sub-analytic
Lagrangian cycles on $T^*X$ with integral coefficients. (One way of
seeing this is to use the fact that any stratification can be refined
to a triangulation, where monodromy no longer plays a role.)

\begin{remark}
Using Kashiwara's definition it is not immediately obvious
that one obtains a cycle. However, for standard sheaves on simplices
this follows by the limit construction given in \cite{SV1}. As these
objects generate $K(X;\Bbbk)$ it follows that one obtains a cycle in general. 
\end{remark}

Consider the functors of extension of scalars:
\[
D^b_c(X, \FM_p)
\stackrel{(-)\stackrel{L}{\otimes}_\Z\FM_p}{\longleftarrow}
D^b_c(X, \Z) \stackrel{{(-) \otimes_{\Z} \CM}}{\longrightarrow} D^b_c(X, \CM).
\]
They induce maps on the corresponding Grothendieck groups
\[
K(X;\FM_p) \stackrel{f}{\from} K(X;\ZM) \stackrel{c}{\to} K(X;\CM)
\]
which commute with $\CC$. We summarize this discussion by an analogue of Lemma \ref{lem:commute}:

\begin{lem}
 $f$ and $c$ are isomorphisms commuting with $\CC$ and hence all the groups $K(X;\FM_p)$, $K(X;\ZM)$, $K(X;\CM)$, and  $\LC_X$ are canonically isomorphic. 
\end{lem}

(One can see directly that $f$ and $c$ are isomorphisms by again
refining to a triangulation.)

  Finally, let us note that an alternative definition of the characteristic cycle has been given
  by Kashiwara and Schapira \cite[\S 9.4]{KS}
under the assumption that $\Bbbk$ is
  a field of characteristic zero. It has the disadvantage that it does
  not work well with finite characteristic coefficients.

\section{Examples}

\subsection{The flag variety of type $B_2$:} Let $G$ be a simple complex
algebraic group of  type $B_2$, $B \subset G$ a Borel subgroup and $X
= G/B$ its flag variety. Let $W$ be the Weyl group of $G$ with simple
reflections $s$ and $t$, corresponding to the short and long simple
roots respectively. Consider the Bruhat decompositon
\[
X = \bigoplus_{w \in W} BwB/B = \bigsqcup_{w \in W} X_w.
\]
Then the Schubert varieties $\overline{X_w}$ are smooth unless $w =
sts$. In this case
\[
\CC(\Ic(\overline{X_{sts}})) = [\overline{T_{sts}^*X}] + [\overline{T_{s}^*X}].
\]
This is the first example of a reducible characteristic variety of a
Schubert variety, and is due to  Kashiwara and Tanisaki \cite{KT}.

The singularity of $X_{sts}$ along $X_s$ is smoothly equivalent to
a Kleinian surface singularity of type $A_1$. Hence one has (see
\cite[\S 2.5.3]{J})
\[
d^p_{sts, s} = \begin{cases} 1 &\text{if $p =2$,} \\ 0
  &\text{otherwise}. \end{cases}
\]

\subsection{Flag varieties of type $A$:} Let $G = SL_n(\CM)$, $B
\subset G$ a Borel subgroup and $W=S_n$ be the Weyl group. As in the
previous example we consider the Bruhat decomposition
\[
X = \bigoplus_{w \in W} BwB/B = \bigsqcup_{w \in W} X_w.
\]
Using the computer calculations in \cite[Section 5.1]{W} one can show
there are no decomposition numbers if $n < 8$. If $n = 8$ then one can use
\cite[Section 5.2]{W} to conclude that all decomposition numbers are
zero except for 38 possible
exceptions. Building of work of Braden \cite[Appendix A]{W} 
the second author has performed computer calculations to verify
that $\ic(\overline{X_w};\ZM) \otimes_\ZM^L \FM_2$ is not simple in these
remaining 38 cases. Our main theorem now implies that all of these
Schubert varieties have reducible characteristic variety.

One of these 38 Schubert varieties (the Schubert
variety corresponding to the permutation 62845173) has a singularity
smoothly equivalent to a singularity in a quiver variety of type
$A_5$. This singularity was used by Kashiwara and Saito \cite[Theorem
7.2.1]{KS} to provide a
counterexample to a conjecture of Kazhdan and Lusztig on the
irreducibility of characteristic varieties in type $A$ quiver
varieties and Schubert varieties. The second author \cite{W2} has recently
noticed that this example also gives a counterexample to the
Kleshchev-Ram conjecture \cite{KR}.


\end{document}